\documentclass{amsart}

\newtheorem{thm}{Theorem}[section]
\newtheorem{lemma}[thm]{Lemma}
\newtheorem{prop}[thm]{Proposition}
\newtheorem{cor}[thm]{Corollary}
\newtheorem{fact}[thm]{Fact}
\newtheorem{quest}[thm]{Question}
\newtheorem{conj}[thm]{Conjecture}

\theoremstyle{definition}
\newtheorem{defn}[thm]{Definition}
\newtheorem{exam}[thm]{Example}
\newtheorem{discuss}[thm]{Discussion}

\theoremstyle{remark}
\newtheorem{remark}{Remark}

\numberwithin{equation}{section}

\newcommand{\ehk}{e_{\rm HK}}
\newcommand{\shk}{s_{\rm HK}}
\newcommand{\e}{e}
\newcommand{\ord}{\mathrm{ord}\,}
\newcommand{\Assh}{\mathrm{Assh}\,}

\newcommand{\chara}{\mathrm{char}\,}
\newcommand{\degree}{\mathrm{deg}\,}
\newcommand{\emb}{\mathrm{emb}\,}
\newcommand{\height}{\mathrm{height}\,}
\newcommand{\gr}{\mathrm{gr}}
\newcommand{\rank}{\mathrm{rank}\,}
\newcommand{\vol}{\mathrm{vol}\,}

\newcommand{\frm}{{\mathfrak m}}
\newcommand{\fraM}{{\mathfrak M}}
\newcommand{\frn}{{\mathfrak n}}
\newcommand{\frp}{{\mathfrak p}}
\newcommand{\bbF}{\ensuremath{\mathbb F}}
\newcommand{\bbZ}{\ensuremath{\mathbb Z}}

\newcommand{\bbM}{\ensuremath{\mathbb M}}
\newcommand{\bbP}{\ensuremath{\mathbb P}}
\newcommand{\bbR}{\ensuremath{\mathbb R}}
\newcommand{\bbQ}{\ensuremath{\mathbb Q}}
\newcommand{\II}{\ensuremath{\mathcal I}}
\newcommand{\HH}{\ensuremath{\mathcal H}}
\newcommand{\PP}{\ensuremath{\mathcal P}}
\newcommand{\CC}{\ensuremath{\mathcal C}}
\newcommand{\LL}{\ensuremath{\mathcal L}}

\begin{document}

\title{Hilbert-Kunz multiplicity of three-dimensional 
local rings}

\author{Kei-ichi Watanabe}
\address{Department of Mathematics, 
College of Humanities and Sciences,
Nihon University, Setagaya-ku, Tokyo 156--0045, Japan}
\email{watanabe@math.chs.nihon-u.ac.jp}
\thanks{The first author was supported in part by  
Grant aid in Scientific
Researches, \# 13440015 and \# 13874006.}  
\author{Ken-ichi Yoshida}
\address{Graduate School of Mathematics, 
Nagoya University,
Chikusa-ku, Nagoya  464--8602, Japan}
\email{yoshida@math.nagoya-u.ac.jp}
\thanks{The second author was supported in part by 
NSF Grant \#14540020.}
\subjclass{Primary 13D40, 13A35; 
 Secondary, 13H05, 13H10, 13H15}




\begin{abstract}
In this paper, we investigate a lower bound (say $\shk(p,d)$)
on Hilbert-Kunz multiplicities
for non-regular unmixed local rings of Krull dimension $d$
with characteristic $p>0$.
Especially, we focus three-dimensional local rings.
In fact, as a main result, we will prove that
$\shk(p,3) = 4/3$ and that a three-dimensional complete
local ring of Hilbert-Kunz multiplicity $4/3$ is isomorphic
to the non-degnerate quadric hyperplanes 
$k[[X,Y,Z,W]]/(X^2+Y^2+Z^2+W^2)$
under mild conditions.
\par
Furthermore, we pose a generalization of the main theorem
to the case of $\dim A \ge 4$ as a conjecture, and show that
it is also true in case of $\dim A = 4$ using the similar 
method as in the proof of the main theorem.
\end{abstract}

\maketitle

\section*{Introduction}

Let $A$ be a commutative Noetherian ring of 
characteristic $p>0$ with unity.
In \cite{Ku1}, Kunz proved the following theorem, 
which gives a characterization of
regular local rings of positive characteristic.

\vspace{2mm}\par\noindent
{\bf Kunz' Theorem (\cite{Ku1}).}
Let $(A,\frm,k)$ be a local ring of characteristic $p>0$.
Then the following conditions are equivalent$:$
\begin{enumerate}
 \item $A$ is a regular local ring.
 \item $A$ is reduced and is flat over 
       a subring $A^p=\{a^p\,|\, a\in A\}$.
  In other words, the Frobenus map 
  $F:A\to A\;(a \mapsto a^{p})$ is flat.
 \item $l_A(A/\frm^{[q]}) = q^d$ for any $q=p^e, e \ge 1$,
 where $\frm^{[q]} = (a^q \,|\, a\in \frm)$ 
 and $l_A(M)$ denotes the length of an $A$-module $M$.
\end{enumerate}

\par 
Furthermore, in \cite{Ku2}, Kunz observed that
$l_A(A/\frm^{[q]})/q^d$ $(q=p^e)$
is a reasonable measure for singularity of a local ring.
Based on the idea of Kunz, Monsky \cite{Mo1} proved that
there exists a constant $c=c(A)$ such that
\[
l_A(A/\frm^{[q]}) = c q^d +O(q^{d-1})
\]
and defined the notion of 
{\it Hilbert-Kunz multiplicity} by $\ehk(A)=c$.
In 1990's, Han and Monsky \cite{HM} have given  
an algorism to compute the Hilbert-Kunz multiplicity 
for any hypersurface of Briskorn--Fermat type
\[
 A=k[X_0,\ldots,X_n]/(X_0^{d_0}+\cdots + X_n^{d_n}).
\]
See e.g. \cite{BC,BCP,Co,WY3} about the other examples.
Hochster and Huneke \cite{HH1}
have given \lq\lq Length Criterion for Tight Closure''
in terms of Hilbert-Kunz multiplcity 
(see Theorem \ref{LCforTC})
and indicated the close relation between 
the tight closure and the Hilbert-Kunz multiplicity.
In \cite{WY1}, the authors have proved the theorem 
which gives a characterization of regular local rings
in terms of Hilbert-Kunz multiplicity as follows$:$

\vspace{2mm} \par\noindent
{\bf Theorem A (\cite[Theorem (1.5)]{WY1}).}
Let $(A,\frm,k)$ be an unmixed local ring of 
positive characteristic.
Then $A$ is regular if and only if $\ehk(A) =1$.

\par \vspace{2mm}
Many researchers have tried to improve this theorem. 
For example, Blickle and Enescu \cite{BE} recently 
proved the following theorem$:$

\vspace{2mm}\par\noindent
{\bf Theorem B (\cite{BE}).}
Let $(A,\frm,k)$ be an unmixed local ring of 
characteristic $p>0$.
Then the following statements hold$:$
\begin{enumerate}
 \item If $\ehk(A) < 1 + \frac{1}{d!}$, 
  then $A$ is Cohen--Macaulay and F-rational.
 \item If $\ehk(A) < 1 + \frac{1}{p^d d!}$, 
  then $A$ is regular.
\end{enumerate}

\par \vspace{2mm}
So it is natural to consider the following problem.

\vspace{2mm} \par\noindent
{\bf Problem C.}
Let $d \ge 2$ be any integer.
Determine the value
\[
\shk(p,d):=\inf \{\ehk(A) \,|\; 
\mbox{$A$ is a non-regular unmixed local ring of
$\chara A =p$} \}.
\]
Also, characterize the local rings $A$ 
for which $\ehk(A) =  \shk(p,d)$ hold.

\par \vspace{2mm}
In case of one-dimensional local rings, 
it is easy to answer to this
problem.
In fact, $\shk(p,1) =2$; 
$\ehk(A) = 2$ if and only if $e(A) =2$.
In case of two-dimensional Cohen--Macaulay local rings,
the authors \cite{WY2} have given a complete answer to 
this problem.
Namely, we have $\shk(p,2) = \frac{3}{2}$
 by the theorem below.

\vspace{2mm} \par\noindent
{\bf Theorem D.} (See also Corollary \ref{2dimHK}.)
Let $(A,\frm,k)$ be a two-dimensional Cohen--Macaulay 
local ring of positive characteristic.
Put $e=e(A)$, the multiplicity of $A$.
Then the following statements hold$:$
\begin{enumerate}
 \item $\ehk(A) \ge \frac{e+1}{2}$.
 \item Suppose that $k = \overline{k}$.
  Then $\ehk(A) = \frac{e+1}{2}$ holds 
  if and only if the associated graded ring 
  $\gr_{\frm}(A)$ is isomorphic to the Veronese 
  subring $k[X,Y]^{(e)}$.
\end{enumerate}

\par \vspace{2mm}
In the following, let us explain the 
organization of this paper.
In Section 1, we recall the notions of 
Hilbert-Kunz multiplicity
and tight closure etc.~ and gather several 
fundamental properties of them.
In particular, Goto--Nakamura's theorem 
(Theorem \ref{GotoNak}) is
important because it plays a central role 
in the proof of main result (Theorem \ref{Main3}) 
and the others.

\par
In Section 2, we give a key result to estimate 
Hilbert-Kunz multiplicities for local rings 
in lower dimension.
Indeed, Theorem \ref{Key} is a refinement of the argument
in \cite[Section 2]{WY2}.
Also,  the point of our proof is to estimate
$l_A(\frm^{[q]}/J^{[q]})$
(where $J$ is a minimal reduction of $\frm$) using
volumes in $\bbR^d$.
\par
In Section 3, we prove the following theorem
as the main result in this paper.  

\vspace{2mm}\par \noindent
{\bf Theorem \ref{Main3}.}
Let $(A,\frm,k)$ be a three-dimensional unmixed 
local ring of characteristic $p >0$.
Then the following statements hold.
\begin{enumerate}
 \item If $A$ is not regular, 
  then $\ehk(A) \ge \frac{4}{~3~}$.
 \item Suppose that $k = \overline{k}$ and 
  $\chara k \ne 2$.
  Then the following conditions are equivalent$:$
 \begin{enumerate}
  \item $\ehk(A) = \frac{4}{~3~}$.
  \item $\widehat{A} \cong k[[X,Y,Z,W]]/(X^2+Y^2+Z^2+W^2)$.
 \end{enumerate}
\end{enumerate}

\par \noindent
Also, we study lower bounds on $\ehk(A)$ for local 
rings $A$ having a given (small) multiplicity $e$.
In particular, we will prove that any three-dimensional 
unmixed local ring $A$ with $\ehk(A) < 2$ is F-rational. 

\par
In Section 4, we consider a generalization 
of Theorem \ref{Main3} and pose the following conjecture.

{\bf Conjecure \ref{Conj-Gen}.}
Let $d \ge 1$ be an integer and $p > 2$ a prime number.
Put
\[
  A_{p,d} := \overline{\bbF_p}[[x_0,x_1,\ldots,x_d]]/
 (x_0^2+\cdots + x_d^2).
\]
Let $(A,\frm,k)$ be a $d$-dimensional unmixed 
local ring with $k =\overline{\bbF_p}$.
Then the following statements hold.
\begin{enumerate}
 \item If $A$ is not regular, 
 then $\ehk(A) \ge \ehk(A_{p,d}) \ge 1+\frac{c_d}{d!}\;$ 
 (See \ref{Conj-Gen} for the definition of $c_d$). 
  In particular, $\shk(p,d) = \ehk(A_{p,d})$.
 \item If $\ehk(A) = \ehk(A_{p,d})$, then 
  the $\frm$-adic completion $\widehat{A}$ of $A$ is isomorphic 
  to $A_{p,d}$ as local rings. 
\end{enumerate}

\par \vspace{2mm}
Also, we prove that this is true in case of $\dim A =4$.
Namely we will prove the following theorem.

\par \vspace{2mm}\par \noindent
{\bf Theorem \ref{4Main}.}
Let $(A,\frm,k)$ be a four-dimensional 
unmixed local ring of characteristic $p >0$. 
Also, suppose that $k=\overline{k}$ and $\chara k \ne 2$. 
Then  $\ehk(A) \ge \frac{5}{~4~}$ if $e(A) \ge 3$, 
where $e(A)$ denotes the multiplicity of $A$.
Also, the following statements hold.
\begin{enumerate}
 \item If $A$ is not regular, then $\ehk(A) \ge
 \ehk(A_{p,4}) = \frac{29p^2+15}{24p^2+12}$.
 \item The following conditions are equivalent$:$
\begin{enumerate}
   \item Equality holds in $(1)$.
   \item $\ehk(A) < \frac{~5~}{4}$.
   \item $\widehat{A}$ is isomorphic to $A_{p,4}$. 
\end{enumerate}
\end{enumerate}

\section{Preliminaries} 

\par
Throughout this paper, let $A$ be a commutative Noetherian 
ring with unity.
Furthermore, we assume that $A$ has a positive 
characteristic $p$, that is,
it contains a prime field $\bbF_p= \bbZ/p\bbZ$, 
unless specified.
For every positive integer $e$, let $q = p^e$.
If $I$ is an ideal of $A$, then 
$I^{[q]} = (a^q \,|\, a \in I)A$.
Also, we fix the following notation$:$
$l_A(M)$ (resp. $\mu_A(M))$ denotes the length
(resp. the minimal number of generators) of $M$
for any finitely generated $A$-module $M$.

\par
First, we recall the notion of Hilbert-Kunz multiplicity
(see \cite{Ku1,Ku2,Mo1}), which plays a central
role in this paper.
Also, see \cite{Ma} or \cite{Na} for usual multiplicity.

\begin{defn}[Multiplicity, Hilbert-Kunz Multiplicity]
Let $(A,\frm,k)$ be a local ring of characteristic $p>0$ 
with $\dim A = d$.
Let $I$ be an $\frm$-primary ideal of $A$, 
and let $M$ be a finitely generated $A$-module.
The {\it $($Hilbert-Samuel$)$ multiplicity\/} $\e(I,M)$ 
of $I$ with respect to $M$ is defined by
\[
  \e(I,M) = \lim_{n \to \infty} 
  \frac{d!}{n^d} \;l_A(M/I^n M)
\]
The {\it Hilbert-Kunz multiplicity\/} $\ehk(I,M)$ of $I$ 
with respect to $M$ is defined by
\[
  \ehk(I,M) = \lim_{q \to \infty} 
  \frac{l_A(M/I^{[q]} M)}{q^d}.
\]
By definition, we put $\e(I) = \e(I,A)$ 
(resp. $\ehk(I) = \ehk(I,A)$)
and $\e(A) = \e(\frm)$ 
(resp. $\ehk(A) = \ehk(\frm)$).
\end{defn}

\par
We recall several basic results on 
Hilbert-Kunz multiplicity.

\begin{prop}[Fundamental Properties 
{\rm (cf. \cite{Hu,Ku1,Ku2,Mo1,WY1})}]
\label{Fund}
Let $(A,\frm,k)$ be a local ring of characteristic $p >0$.
Let $I$, $I'$ be $\frm$-primary ideals of $A$, and 
let $M$ be a finitely generated $A$-module.
Then the following statements hold.
\begin{enumerate}
 \item If $I \subseteq I'$, then $\ehk(I) \ge \ehk(I')$.
 \item $\ehk(A) \ge 1$.
 \item $\dim M < d$ if and only if $\ehk(I,M) = 0$.
 \item If $0 \to L \to M \to N \to 0$ is a short exact 
  sequence of finitely generated $A$-modules, then
\[
  \ehk(I,M) = \ehk(I,L) + \ehk(I,N).
\]
 \item {\rm (Associative Formula)}
\[
  \ehk(I,M) = \sum_{\frp \in \Assh(A)} \ehk(I,A/\frp)
\cdot l_{A_\frp}(M_{\frp}),
\]
where $\Assh(A)$ denotes the set of minimal prime 
ideals $\frp$ of $A$ with $\dim A/\frp = \dim A$.
 \item If $J$ is a parameter ideal of $A$, 
  then $\ehk(J) = \e(J)$.
  In particular, if $J$ is a minimal reduction of $I$
  $($i.e., $J$ is a paramater ideal which is 
  contained in $I$ and $I^{r+1} = JI^r$ for some 
  integer $r \ge 0$$)$, then $\ehk(J) = \e(I)$.
 \item If $A$ is regular, then $\ehk(I) = l_A(A/I)$.
 \item {\rm (Localization)}
  $\ehk(A_{\frp}) \le \ehk(A)$ holds for
  any prime ideal $\frp$ such that 
  $\dim A/\frp + \height \frp = \dim A$.
 \item If $x \in I$ is $A$-regular, 
  then $\ehk(I) \le \ehk(I/xA)$.
 \item If $(A,\frm) \to (B,\frn)$ is a flat local ring 
  homomorphism such that $B/\frm B$ is a field, 
  then $\ehk(I) = \ehk(IB)$.
\end{enumerate}
\end{prop}

\begin{remark}
 Also, the similar result as above (except (6),(7)) holds for
usual multiplicities.
\end{remark}

\par
Let $(A,\frm,k)$ be any local ring of positive dimension. 
{\it The associated graded ring} $\gr_{\frm}(A)$ of $A$ 
with respect to $\frm$ is defined as follows$:$
\[
\gr_{\frm}(A) :=  \bigoplus_{n \ge 0} \frm^n/\frm^{n+1}.
\]
Then $G=\gr_{\frm}(A)$ is a homogeneous $k$-algebra such that
$\fraM := G_{+}$ is
the unique homogeneous maximal ideal of $G$.
If $\chara A = p >0$ and $\dim A =d$, then
$G_{\fraM}$ is also a local ring of characteristic
$p$ with $\dim G_{\fraM} = d$.

\begin{prop} {\rm (\cite[Theorem (2.15)]{WY1}, \cite{HuY})} 
\label{grHK}
Let $(A,\frm,k)$ be a local ring of characteristic $p >0$.
Let $G = \gr_{\frm}(A)$ the associated graded ring of $A$ 
with respect $\frm$ as above. 
Then $\ehk(A) \le \ehk(G_{\fraM}) \le e(A)$.
\end{prop}

\begin{remark}
We use the same notation as in the above proposition.
Although $e(A) = e(G_{\fraM})$ always holds,
$\ehk(A) = \ehk(G_{\fraM})$ seldom holds.
\end{remark}

\begin{prop}[{\rm cf. \cite{Hu}}]   \label{IneqMul}
Let $(A,\frm,k)$ be a local ring of characteristic $p >0$
with $d= \dim A$.
Let $I$ be an $\frm$-primary ideal of $A$.
Then
\[
 \frac{\e(I)}{d!} \le \ehk(I) \le \e(I).
\]
Also, if $d \ge 2$, then the inequality 
in the left-hand side is strict$;$ see \cite{Ha}.
\end{prop}

\par
We say that a local ring $A$ is {\it unmixed\/} if
$\dim \widehat{A}/\frp = \dim \widehat{A}$ holds for
any associated prime ideal $\frp$ of $\widehat{A}$.
The following theorem is an analogy of 
Nagata's theorem; see \cite[(40.6)]{Na}.
Furthermore, it is a starting point in this article.

\begin{thm}[{\rm \cite[Theorem (1.5)]{WY1}}]  
\label{Regular}
Let $(A,\frm,k)$ be an unmixed local ring of positive
characteristic. 
Then $A$ is regular if and only if $\ehk(A) = 1$.
\end{thm}

\par
It is not so easy to compute Hilbert-Kunz 
multiplicities in general.
However, one has simple formulas for them in 
case of quotient singularities and in case of 
binomial hypersurfaces; see below or \cite[Theorem 3.1]{Co}. 

\begin{thm}[{\rm cf. \cite[Theorem (2.7)]{WY1}}] 
\label{Quot}
Let $(A,\frm) \hookrightarrow (B,\frn)$ be a 
module-finite extension of local domains.
Then for every $\frm$-primary ideal $I$ of $A$, 
we have
\[
  \ehk(I) = \frac{\ehk(IB)}{[Q(B):Q(A)]} 
  \cdot [B/\frn : A/\frm],
\]
where $Q(A)$ denotes the fraction field of $A$.
\end{thm}

\par
Now let us see some examples of Hilbert-Kunz 
multiplicities which are given by 
the above formula.
First, we consider the Veronese 
subring $A$ defined by
\[
  A = k[[X_1^{i_1}\cdots X_d^{i_d}
\,|\,i_1,\ldots,i_d \ge 0, \sum i_j = r]].
\]
Applying Theorem \ref{Quot} to 
$A \hookrightarrow B = k[[x,y]]$, we get
\begin{equation}
  \ehk(A) = \frac{1}{~r~} \binom{d+r-1}{r}.
\end{equation}
In particular, if $d =2$, $r =e(A)$, 
then $\ehk(A) = \frac{e(A)+1}{2}$.

\par
Next, we consider the homogeneous coordinate ring of
quadric hyperplanes in $\bbP^3_k$.
Let $k$ be a field of characteristic $p > 2$, and let
$R$ be the homogeneous coordinate ring of
the hyperquadric $Q$ defined by $q=q(X,Y,Z,W)$.
Put $\fraM = R_{+}$, the unique homogeneous 
maximal ideal of $R$, and
$A = R_{\fraM} \otimes_k \overline{k}$.
By suitable coordinate transformation, 
we may assume that
$A$ is isomorphic to one of the following rings$:$
\begin{equation}
\begin{cases}
k[[X,Y,Z,W]]/(X^2),    & \mbox{if}\;\; \rank(q) =1, \\
k[[X,Y,Z,W]]/(X^2-YZ), & \mbox{if}\;\; \rank(q) =2, \\
k[[X,Y,Z,W]]/(XY-ZW),  & \mbox{if}\;\; \rank(q) =3. \\
\end{cases}
\end{equation}
Then $\ehk(A) = 2$, $\frac{3}{~2~}$, 
or $\frac{4}{~3~}$, respectively.

\par
In order to state other important properties of 
Hilbert-Kunz multiplicities,
the notion of tight closure is very important. 
See \cite{HH1,HH2, Hu} for definition and the 
fundamental properties of tight closure. 
In particular, the notion  of $F$-rational ring 
is essential in our argument.  

\begin{defn}[{\rm \cite{FW,HH1,HH2}}] 
\label{Fregdefn}
Let $(A,\frm,k)$ be a local ring of 
characteristic $p >0$.
We say that $A$ is {\it weakly F-regular} 
(resp. {\it F-rational\/})
if every ideal (resp. every paramater ideal) 
is tightly closed.
Also, $A$ is {\it F-regular} 
(resp. {\it F-rational}) if
any local ring of $A$ is weakly F-regular 
(resp. F-rational).
\end{defn}

\par 
Note that an $F$-rational local ring is normal and  
Cohen-Macaulay. 

\par
Hochster and Huneke have given the following 
criterion of tight closure 
in terms of Hilbert-Kunz multiplicity.

\begin{thm}[Length Criterion of Tight Closure
{\rm (cf. \cite[Theorem 8.17]{HH1})}]
\label{LCforTC}
Let $I \subseteq J$ be $\frm$-primary ideals of 
a local ring $(A,\frm,k)$ of characteristic $p >0$.
\begin{enumerate}
 \item If $I^{*} = J^{*}$, then $\ehk(I) = \ehk(J)$.
 \item Suppose that $A$ is excellent, reduced 
  and equidimensional.
  Then the converse of $(1)$ is also true.
\end{enumerate}
\end{thm}

\par
The following theorem plays an important role 
in studying Hilbert-Kunz multiplicities for 
non-Cohen--Macaulay local rings.

\begin{thm}[Goto--Nakamura {\rm \cite{GN}}] 
\label{GotoNak}
Let $(A,\frm,k)$ be an equidimensional local ring
which is a homomorphic image of a 
Cohen--Macaulay local ring 
of characteristic $p >0$.
Then
\begin{enumerate}
 \item If $J$ is a parameter ideal of $A$, 
  then $e(J) \ge l_A(A/J^{*})$.
 \item Suppose that $A$ is unmixed. 
 If $e(J) = l_A(A/J^{*})$, then
 $A$ is F-rational $($and thus$)$ is 
 Cohen--Macaulay.
\end{enumerate}
\end{thm}

\par
The next lemma is well-known in case of
Cohen--Macaulay local rings (e.g. see \cite{Hu}).

\begin{cor} \label{Mul2}
Let $(A,\frm,k)$ be an unmixed local ring of 
characteristic $p>0$.
Suppose that $e(A) =2$.
Then $\widehat{A}$ is F-rational 
if and only if $\ehk(A) < 2$.
When this is the case, $A$ is a hypersurface.
\end{cor}

\begin{proof}[\quad Proof]
Since any Cohen--Macaulay local ring of 
multiplicity $2$ is a hypersurface, 
it suffices to prove the first statement. 
\par 
We may assume that $A$ is complete and $k$ is infinite. 
We can take a minimal reduction $J$ of $\frm$.
First, suppose that $\ehk(A) <2$.
Then we show that $A$ is Cohen--Macaulay, F-rational.
By Goto--Nakamura's theorem,
we have $2 = e(J) \ge l_A(A/J^{*})$.
If equality does not hold, $l_A(A/J^{*}) =1$, 
that is, $J^{*} = \frm$.
Then $\ehk(A) = \ehk(J^{*}) = \ehk(J) = e(J) =2$ 
by Proposition \ref{Fund}.
This is a contradiction.
Hence $e(J) = l_A(A/J^{*})$.
By Goto--Nakamura's theorem again, 
we obtain that $A$ is Cohen--Macaulay, F-rational.
\par
Conversely, suppose that $A$ is 
complete F-rational.
Then since $A$ is Cohen--Macaulay 
and $J^{*} = J \ne \frm$,
we have $\ehk(A) < \ehk(J) =e(J) =2$ 
by Length Criterion for Tight Closure.
\end{proof}

\par
The next question is open in general.
However, we will show that it is true for 
$\dim A \le 3$; see Section 2.

\begin{quest} \label{Ques2}
If $A$ is an unmixed local ring with $\ehk(A) < 2$, 
then is it F-rational?
\end{quest}

\section{Estimate of Hilbert-Kunz multiplicity}

\par
In this section, we will prove the key result to find a 
lower bound of Hilbert-Kunz multiplicities.
Actually, it is a refinement of the argument 
which appeared in \cite[Section 5]{WY1} 
or in \cite[Section 2]{WY2}.
The point is to use the tight closure 
$J^{*}$ instead of \lq\lq a parameter ideal $J$ itself''.
This enables us to investigate Hilbert-Kunz 
multiplicities of non-Cohen--Macaulay local rings.
In Sections $3$, $4$, we will
apply our method to unmixed local rings 
with $\dim A =3,\,4$.

\par
Before stating our theorem, 
we introduce the following notation$:$
For any positive real number $s$, we put
\[
 v_s := \vol \left\{(x_1,\ldots,x_d) \in [0,1]^d 
 \,\bigg|\, \sum_{i=1}^d x_i \le s\right\}, \quad
 v_s':=1-v_s,
\]
where $\vol(W)$ denotes the volume of 
$W \subseteq \bbR^d$.
Then it is easy to see the following fact.

\begin{fact} \label{VolF}
Let $d \ge 1$ be an integer, and let $s$ be a
 positive real  number.
Using the same notation as above, we have
\begin{enumerate}
 \item $v_s + v_s' = 1$.
 \item $v_{d-s}' = v_s$.
 \item $v_{d/2} = v_{d/2}' = \frac{1}{~2~}$.
 \item If $0 \le s \le 1$, 
then $v_s = \frac{s^d}{d!}$.
\end{enumerate}
\end{fact}

\par
Using the above notaion, a key result in this paper 
can be written as follows$:$

\begin{thm} \label{Key}
Let $(A,\frm,k)$ be an unmixed local ring of 
characteristic $p>0$.
Put $d = \dim A \ge 1$.
Let $J$ be a minimal reduction of $\frm$, and let
$r$ be an integer with $r \ge \mu_A(\frm/J^{*})$, 
where $J^{*}$ denotes the tight closure of $J$.
Also, let $s \ge 1$ be a rational number.
Then we have
\begin{equation} \label{KeyEq}
  \ehk(A) \ge e(A)
 \left\{v_s - r \cdot \frac{(s-1)^d}{d!}\right\}.
\end{equation}
\end{thm}

\begin{remark}
When $1 \le s \le 2$, the right-hand side in 
Equation (\ref{KeyEq}) is equal to 
$e(A)(v_s - r \cdot v_{s-1})$.
\end{remark}

\par
Before proving the theorem, we need the following lemma.  
In what follows, for any positive real number $\alpha$,
we define $I^{\alpha} := I^{n}$, where $n$ is 
the minimum integer which does not exceed $\alpha$.

\begin{lemma} \label{Volume}
Let $(A,\frm,k)$ be an unmixed local ring
of characteristic $p>0$ with $\dim A = d \ge 1$.
Let $J$ be a parameter ideal of $A$.
Using the same notation as above, we have
\[
 \lim_{q \to \infty} \frac{l_A(A/J^{sq})}{q^d} 
= \frac{e(J) s^d}{d!},
 \qquad
 \lim_{q\to \infty} l_A
 \left(\frac{J^{sq}+J^{[q]}}{J^{[q]}}\right)
 = e(J) \cdot v_s'.
\]
\end{lemma}

\begin{proof}[\quad Proof]
We may assume that $A$ is complete.
Let $x_1,\ldots,x_d$ be a system of parameters 
which generates $J$,
and put $R := k[[x_1,\ldots,x_d]]$, 
$\frn = (x_1,\ldots,x_d)R$.
Then $R$ is a complete regular local ring and
$A$ is a finitely generated $R$-module 
with $A/\frm = R/\frn$.
Since the assertion is clear in case of 
regular local rings,
it suffices to show the following claim.

\begin{description}
 \item[Claim] Let $\II = \{I_q\}_{q=p^e}$ be 
a set of ideals of $A$ which
satisfies the following conditions$:$
\begin{enumerate}
\item For each $q=p^e$, $I_q = J_q A$ holds 
for some ideal $J_q \subseteq
R$.
\item There exists a positive integer $t$ such that
$\frn^{tq} \subseteq J_q$ for all $q=p^e$.
\item $\lim_{q\to \infty} l_R(R/J_q)/q^d$ exists.
\end{enumerate}
Then
\[
\lim_{q\to \infty} \frac{l_A(A/I_q)}{q^d}
= e(J)\cdot \lim_{q\to \infty} \frac{l_R(R/J_q)}{q^d}.
\]
\end{description}

\par
In fact, since $A$ is unmixed, it is a torsion-free 
$R$-module of rank $e:=e(J)$.
Take a free $R$-module $F$ of rank $e$ such that 
$A_W \cong F_W$, where
$W = A \setminus \{0\}$.
Since $F$ and $A$ are both torsion-free, 
there exist the following short
exact sequences of finitely generated $R$-modules$:$
\[
 0 \to F \to A \to C_1 \to 0,
\qquad  0 \to A \to F \to C_2 \to 0,
\]
where $(C_1)_W = (C_2)_W =0$.
In particular, $\dim C_1 < d$ and $\dim C_2 <d$.
\par
Applying the tensor product $-\otimes_R R/J_q$ 
to the above two exact sequences,
respectively, we get
\begin{eqnarray*}
 l_R(A/I_q) & \le &l_R(F/J_qF) + l_R(C_1/J_q C_1),  \\
 l_R(F/J_q F) & \le & l_R(A/I_q) + l_R(C_2/J_q C_2).
\end{eqnarray*}
In general, if $\dim_R C < d$, then
\[
\frac{l_R(C/J_q C)}{q^d} \le 
\frac{l_R(C/\frn^{tq} C)}{q^d} \to 0 \quad
(q\to \infty).
\]
Thus the required assertion easily follows 
from the above observation.
\end{proof}

\begin{proof}[\quad {\bf Proof of Theorem \ref{Key}}]
For simplicity, we put $L = J^{*}$.
We will give an upper bound 
of $l_A(\frm^{[q]}/J^{[q]})$. 
First, we have the following inequality$:$ 
\begin{eqnarray*}
  l_A(\frm^{[q]}/J^{[q]})
 & \le &l_A\left(\frac{\frm^{[q]}
  +\frm^{sq}}{J^{[q]}}\right) \\
 & = & l_A\left(
  \frac{\frm^{[q]}+\frm^{sq}}{L^{[q]}+\frm^{sq}}\right)
   + \,l_A\left(
   \frac{L^{[q]}+\frm^{sq}}{L^{[q]}+J^{sq}}\right) \\
 & \phantom{=} & + \,l_A
  \left(\frac{L^{[q]}+J^{sq}}{J^{[q]}+J^{sq}}\right)
   +\,
  l_A\left(\frac{J^{[q]}+J^{sq}}{J^{[q]}}\right)
  =:\ell_1+\ell_2+\ell_3+\ell_4.
\end{eqnarray*}
\par
Next, we see that $\ell_1 \le r \cdot 
l_A(A/J^{(s-1)q}) + O(q^{d-1})$.
By our assumption, we can write as 
$\frm= L+ Aa_1+\cdots +Aa_r$.
Since $\frm^{(s-1)q}a_i^q \subseteq 
\frm^{sq} \subseteq \frm^{sq}+L^{[q]}$, 
we have
\begin{eqnarray*}
  \ell_1 = l_A\left(
  \frac{\frm^{[q]}+\frm^{sq}}{L^{[q]}+\frm^{sq}}\right)
 & \le & \sum_{i=1}^rl_A\left(
  \frac{Aa_i^q+L^{[q]}+\frm^{sq}}{L^{[q]}+\frm^{sq}}
  \right) \\
 & = & \sum_{i=1}^r 
  l_A\left(A/(L^{[q]}+\frm^{sq}):a_i^q \right) \\
 & \le & r \cdot l_A(A/\frm^{(s-1)q}).
\end{eqnarray*}
Since $J$ is a minimal reduction of $\frm$,
we have $l_A(\frm^{tq}/J^{tq}) = O(q^{d-1})$.
Thus we have the required inequality.
Similarly, we get
\[
 \ell_2 = l_A\left(
 \frac{L^{[q]}+\frm^{sq}}{L^{[q]}+J^{sq}}\right)
 \le l_A(\frm^{sq}/J^{sq}) = O(q^{d-1}).
\]

\par
Also, we have
$l_A(L^{[q]}/J^{[q]}) = O(q^{d-1})$ by
Length Criterion for Tight Closure.
Hence
\[
 l_A(\frm^{[q]}/J^{[q]}) \le 
r \cdot l_A(A/J^{(s-1)q})
 +l_A\left(
\frac{J^{[q]}+J^{sq}}{J^{[q]}} \right) +O(q^{d-1}).
\]
It follows from the above argument that
\begin{eqnarray*}
 \ehk(J) - \ehk(\frm)
 &\le & r \cdot \lim_{q \to \infty} 
  \frac{l_A(A/J^{(s-1)q})}{q^d}
 + \lim_{q \to \infty} \frac{1}{q^d} 
 l_A\left(\frac{J^{[q]}+J^{sq}}{J^{[q]}}\right) \\
 &=& r \cdot e \cdot \frac{(s-1)^d}{d!} + e \cdot v_s'.
\end{eqnarray*}
Since $\ehk(J)=e(J) = e$, $\ehk(A) = \ehk(\frm)$ 
and $v_s' = 1-v_s$,
we get the required inequality.
\end{proof}

\par
The following fact is known, which gives a lower bound on 
Hilbert-Kunz multiplicities for hypersurface local rings. 

\begin{fact}[{\rm cf. \cite{BC,BCP,WY1}}]  \label{Hyp}
Let $(A,\frm,k)$ be a hypersurface local ring of
characteristic $p>0$ with $d = \dim A \ge 1$.
Then
\[
   \ehk(A) \ge \beta_{d+1} \cdot e(A), 
\]
where $\beta_{d+1}$ is given by the following equivalent formulas$:$
\begin{eqnarray*}
 (a) && \frac{1}{\pi} 
 \int_{-\infty}^{\infty}
 \left(\frac{\sin \theta}{\theta}\right)^{d+1} 
 \,d \theta;  \\
 (b) && \frac{1}{2^d d!} 
 \sum_{\ell = 0}^{[\frac{d}{2}]} (-1)^{\ell}
 (d+1-2\ell)^d
 \binom{d+1}{\ell}; \\
 (c) && \vol \left\{\underline{x} 
 \in [0,1]^d \,\bigg|\, \frac{d-1}{2} \le
 \sum x_i
 \le \frac{d+1}{2}\right\}
 =1 - v_{\frac{d-1}{2}} - v_{\frac{d+1}{2}}'.
\end{eqnarray*}

\begin{table}[ht]
\caption{}\label{eqtable}
\renewcommand\arraystretch{1.5}
\noindent\[
\begin{array}{|c||c|c|c|c|c|c|c|} \hline
d & 0 & 1 & 2 & 3 & 4 & 5 & 6 \\ \hline
\beta_{d+1} & 1 & 1 & \frac{~3~}{4} 
& \frac{2}{~3~} & \frac{115}{192} &
\frac{11}{20} &
\frac{5633}{11520} \\ \hline
\end{array}
\]
\end{table}
\end{fact}

\begin{remark} 
The above inequality is {\it not} best possible 
in general. 
In case of $d \ge 4$, one cannot prove the formula in the above fact  
as a corollary of our theorem. 
See also Proposition \ref{Hyp2-Est} and Theorem \ref{4Main}. 
\end{remark}

\par
Using Theorem \ref{Key} and 
Lemma \ref{Sally-type} of the next section, 
one can prove the following corollary,
which has been already proved in \cite{WY2} 
in the case of Cohen--Macaulay local rings.

\begin{cor}[{\rm cf. \cite{WY2}}] \label{2dimHK}
Let $(A,\frm,k)$ be a two-dimensional unmixed local ring of 
characeristic $p>0$.
Put $e=e(A)$.
Then
\begin{equation} \label{2dimIneq}
 \ehk(A) \ge \frac{e+1}{2}.
\end{equation}
\par
Also, suppose $k=\overline{k}$.
Then equality holds if and only if 
$\gr_{\frm}(A)$ is isomorphic to the Veronese subring 
$k[X,Y]^{(e)} = k[X^e,X^{e-1}Y,\ldots,XY^{e-1},Y^e]$. 
\par
Moreover, if $A$ is not $F$-rational, then we have 
\[
  \ehk(A) \ge \frac{e^2}{2(e-1)}. 
\]
\end{cor}

\begin{exam}[{\rm \cite[Corollary 3.19]{FT}}]
Let $E$ be an elliptic curve over a field $k=\overline{k}$
of characteristic $p>0$, 
and let $\LL$ be a very ample line bundle on $E$ 
of degree $e \ge 3$. 
Let $R$ be the homogeneous coodinate ring 
(the section ring of $\LL$) defined by 
\[
 R = \bigoplus_{n \ge 0} H^0(E,\LL^{\otimes n}). 
\]
Also, put $A = R_{\fraM}$, where $\fraM$ be the unique 
homogeneous maximal ideal of $R$. 
Then we have $\ehk(A) = \frac{e^2}{2(e-1)}$. 
\end{exam}

\vspace{3mm}
\section{Lower bounds in the case of three-dimensional 
local rings}

\par
In this section, we prove the following main theorem 
in this paper, which gives a lower bound
on Hilbert-Kunz multiplicities for non-regular 
unmixed local rings of dimension $3$.

\begin{thm} \label{Main3}
Let $(A,\frm,k)$ be a three-dimensional unmixed local ring of
characteristic $p >0$.
Then
\begin{enumerate}
 \item If $A$ is not regular, 
  then $\ehk(A) \ge \frac{4}{~3~}$.
 \item Suppose that $k = \overline{k}$ and 
  $\chara k \ne 2$.
  Then the following conditions are equivalent$:$
 \begin{enumerate}
  \item $\ehk(A) = \frac{4}{~3~}$.
  \item $\widehat{A} \cong k[[X,Y,Z,W]]/(X^2+Y^2+Z^2+W^2)$.
  \item $\gr_{\frm}(A) \cong k[X,Y,Z,W]/(X^2+Y^2+Z^2+W^2)$.
  That is, $\gr_{\frm}(A) \cong k[X,Y,Z,W]/(XY-ZW)$.
 \end{enumerate}
\end{enumerate}
\end{thm}

\begin{prop} \label{Mul2-3dim}
Let $(A,\frm,k)$ be a three-dimensional unmixed local ring of
characteristic $p >0$. 
If $\ehk(A) < 2$, then $A$ is F-rational.
\end{prop}

\par \vspace{2mm} From now on, we divide the proof of
Theorem \ref{Main3} and Proposition \ref{Mul2-3dim}
into several steps.
The following lemma is an analogy of Sally's theorem$:$
If $A$ is a Cohen--Macaulay local ring, then
$\mu_A(\frm/J) (= \mu_A(\frm) - \dim A) \le e(A)-1$.

\begin{lemma} \label{Sally-type}
Let $(A,\frm,k)$ be an unmixed local ring 
of positive characteristic,
and let $J$ be a minimal reduction of $\frm$.
Then
\begin{enumerate}
 \item $\mu_A(\frm/J^{*}) \le e(A) -1$.
 \item If $A$ is not F-rational, 
  then $\mu_A(\frm/J^{*}) \le e(A) -2$.
\end{enumerate}
\end{lemma}

\begin{proof}[\quad Proof]
We may assume that $A$ is complete and thus
is a homomorphic image of a Cohen--Macalay local ring.
Put $e = e(A)$.
\par
(1) By Theorem \ref{GotoNak},
we have that
$\mu_A(\frm/J^{*}) \le l_A(\frm/J^{*}) \le e(J) -1 = e-1$.
\par
(2) If $A$ is not F-rational, then
$l_A(A/J^{*}) \le e(J) -1 = e-1$.
Thus $\mu_A(\frm/J^{*}) \le e-2$, as required.
\end{proof}

\par \vspace{2mm}
Let $A$ be an unmixed local ring which is not regular.
Put $e = e(A)$, the multiplicity of $A$.
Then $e$ is an integer with $e \ge 2$.
Thus the assertion (1) of Theorem \ref{Main3} follows from
the following lemma.
Also, this implies that if $\ehk(A) = \frac{4}{~3~}$ then $e(A) =2$
without extra assumption.

\begin{lemma} \label{Main3-(1)}
Using the same notation as in Theorem $\ref{Main3}$,
we put $e = e(A)$, the multiplicity of $A$.
\begin{enumerate}
\setlength{\itemsep}{4pt}%
 \item If $e \ge 5$, then $\ehk(A) > 2$
 \item If $e = 4$, then 
       $\ehk(A) \ge \frac{7}{4} > \frac{4}{3}$.
 \item If $e = 3$, then 
       $\ehk(A) \ge \frac{13}{8} > \frac{4}{3}$.
 \item If $e = 2$, then $\ehk(A) \ge \frac{4}{3}$.
\end{enumerate}
\end{lemma}

\begin{remark} 
The lower bounds of $\ehk(A)$ in Lemma \ref{Main3-(1)} are not best possible. 
\end{remark}

\begin{proof}[\quad Proof]
We may assume that $A$ is complete and $k$ is infinite.
By Lemma \ref{Sally-type}(1), 
we can apply Theorem \ref{Key} with $r = e-1$.
Namely, if $1 \le s \le 2$, then
\begin{equation} \label{Ineq-gen}
  \ehk(A) \ge e (v_s - (e-1)v_{s-1})=
  e\left(\frac{s^3}{~6~} - (e+2)\frac{(s-1)^3}{6}\right).
\end{equation}

\par
Define the real-valued function $f_e(s)$ 
by the right-hand side of Eq.~(\ref{Ineq-gen}).
Then one can easily calculate $\max_{1 \le s \le 2} f_e(s)$.
In fact, if $e \ge 2$, then
\[
 \max_{1 \le s \le 2} f_e(s) 
= f\left(\frac{e+2 + \sqrt{e+2}}{e+1}\right)
= \frac{e}{~6~} \left(\frac{e+2+\sqrt{e+2}}{e+1}\right)^2.
\]
But, in order to prove the lemma,
it is enough to use the following values only$:$

\begin{center}

\renewcommand{\arraystretch}{1.5}
\begin{tabular}{|c||c|c|c|} \hline
$s$ & $\frac{~3~}{2}$ & $\frac{~7~}{4}$ & $2$ \\ \hline
$f_e(s)$ & $\frac{e(25-e)}{48}$
& $\frac{e(289-27e)}{384}$ & $\frac{e(6-e)}{6}$ \\ \hline
\end{tabular}
\end{center}

\par \vspace{2mm} \par \noindent
(1) We show that $\ehk(A) >2$ if $e \ge 5$.
If $e \ge 13$, then by Proposition \ref{IneqMul},
\[
   \ehk(A) \ge \frac{~e~}{3!} \ge \frac{13}{~6~}  > 2.
\]
So we may assume that $5 \le e \le 12$.
Applying Eq.~(\ref{Ineq-gen}) for $s=\frac{~3~}{2}$, we get
\[
 \ehk(A) \ge \frac{e(25-e)}{48} \ge \frac{5(25-5)}{48}
 = \frac{25}{12} > 2.
\]

\par \vspace{2mm} \par \noindent
(2) Suppose that $e = 4$.
Actually, applying Eq.~(\ref{Ineq-gen}) 
for $s = \frac{~3~}{2}$, we get
\[
  \ehk(A) \ge \frac{e(25-e)}{48} = \frac{~7~}{4}.
\]
\par \vspace{2mm} \par \noindent
(3) Suppose that $e =3$.
Applying Eq.~(\ref{Ineq-gen}) 
for $s =\frac{~7~}{4}$, we get
\[
  \ehk(A) \ge \frac{e(289-27e)}{384} = \frac{~13~}{8}.
\]

\par \vspace{2mm} \par \noindent
(4) Suppose that $e = 2$.
Applying Eq.~(\ref{Ineq-gen}) for $s =2$,
\[
  \ehk(A) \ge \frac{e(6-e)}{6} = \frac{4}{~3~},
\]
as required.
\end{proof}

\par \vspace{2mm}
Before proving the second assertion of Theorem \ref{Main3},
we prove Proposition \ref{Mul2-3dim}.
For that purpose,
we now focus three-dimensional non-F-rational local rings.
\par
Let $A$ be an unmixed local ring of positive characteristic with
$\dim A =3$.
Put $e = e(A) \ge 2$.
Also, suppose that $A$ is not F-rational.
If $e =2$, then $\ehk(A) =2$ by Lemma \ref{Mul2}.
On the other hand, if $e \ge 5$, then $\ehk(A) >2$ 
by Lemma \ref{Main3-(1)}.
Thus in order to prove Proposition \ref{Mul2-3dim},
it is enough to investigate the cases of $e=3$, $4$.
Namely, Proposition \ref{Mul2-3dim} follows 
from the following lemma.

\begin{lemma} \label{Mul34}
Let $(A,\frm,k)$ be an unmixed local ring of positive 
characteristic with $\dim A =3$. 
Put $e=e(A)$.
Suppose that $A$ is not F-rational.
Then
\begin{enumerate}
 \item If $e=3$, then $\ehk(A) \ge 2$.
 \item If $e=4$, then $\ehk(A) >2$.
\end{enumerate}
\end{lemma}

\begin{proof}[\quad Proof]
We may assume that $k$ is infinite.
Then one can take a minimal reduction (say $J$) of $\frm$.
By Lemma \ref{Sally-type}(2), we can apply 
Theorem \ref{Key} for $r=e-2$. 
Thus if $1 \le s \le 2$, then
\begin{equation} \label{Ineq-nonF}
  \ehk(A) \ge e\left(\frac{s^3}{6} 
  - (e+1)\frac{(s-1)^3}{6}\right).
\end{equation}
\par \noindent
(1) Suppose that $e=3$.
Applying Eq.~(\ref{Ineq-nonF}) for $s =2$, we get
\[
  \ehk(A) \ge \frac{e(7-e)}{6} = 2.
\]
\par \vspace{2mm} \par \noindent
(2) Suppose that $e =4$.
Applying Eq.~(\ref{Ineq-nonF}) for $s =\frac{~7~}{4}$, 
we get
\[
  \ehk(A) \ge \frac{e(316-27e)}{384} = \frac{13}{~6~} >2,
\]
as required.
\end{proof}

\begin{exam} \label{RNS}
Let $R = k[T,xT,xyT,yT,x^{-1}yT,x^{-2}yT,\ldots,x^{-n}yT]$ 
be a rational normal scroll and put
$\frm = (T,xT,xyT,yT,x^{-1}yT,\ldots,x^{-n}yT)$.
Then $A=R_{\frm}$ is a three-dimensional Cohen--Macaulay 
F-rational local domain with $e(A) = n+2$, and
\[
\ehk(A)
 = \frac{~e(A)~}{2} + \frac{e(A)}{6(n+1)}.
\]
\end{exam}

\begin{proof}[\quad Proof]
Let $\PP \subseteq \bbR$ be a convex polytope with vertex set
\[
 \Gamma = \{(0,0),(1,0),(1,1),(0,1),(-1,1),\ldots,(-n,1)\},
\]
and put $\widetilde{\PP} 
:= \big\{(\alpha,1) \in \bbR^3\,|\,\alpha \in
\PP\big\}$
and $d \PP := \big\{d \cdot \alpha \,|\, \alpha \in \PP\big\}$ 
for every integer $d \ge 0$.
Also, if we define a cone $\CC = \CC(\widetilde{\PP})
:= \{r \beta\,|\, \beta \in \widetilde{\PP},\, 
0 \le r \in \bbQ\}$
and regard $R$ as a homogeneous $k$-algebra
with $\degree x = \degree y = 0$ and $\degree T = 1$,
then the basis of $R_d$ corresponds to the set  
$\{(\alpha,d) \in \bbZ^{3}\,|\, 
\alpha \in \bbZ^2 \cap d \PP\} 
= \{(\alpha,d) \in \bbZ^{3}\,|\, 
\alpha \in \bbZ^2\}\cap \CC$.
\setlength{\unitlength}{1mm}
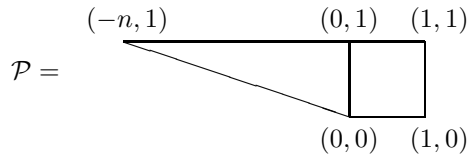
\begin{figure}[ht]
\begin{picture}(100,22)
\put(10,10){\PP \;=}
\put(55,5){\line(0,1){10}} 
\put(65,5){\line(0,1){10}} 
\put(55,5){\line(1,0){10}} 
\put(25,15){\line(1,0){40}} 
\put(25,15){\line(3,-1){30}} 
\put(20,17){$(-n,1)$}
\put(51,17){$(0,1)$}
\put(63,17){$(1,1)$}
\put(51,1){$(0,0)$}
\put(63,1){$(1,0)$}
\end{picture}
\end{figure}
\par
If we put $\Gamma_q 
=\{(0,0),(q,0),(q,q),(0,q),(-q,q),\ldots,(-nq,q)\}$,
then  $\frm^{[q]} = (x^a y^b T^q\,|\,
 (a,b) \in \Gamma_q)$.
Since $[\frm^{[q]}]_d = \sum_{(a,b)\in \Gamma_q} 
R_{d-q} \,x^ay^bT^q$,
we have
\begin{eqnarray*}
 \ehk(A)
 & = & \lim_{q\to \infty} \frac{1}{q^3} 
\,l_A(A/\frm^{[q]}) \\
 & = & \lim_{q\to \infty} \frac{1}{q^3} \,
  \#\bigg\{\bbZ^3 \cap 
\big(\CC \setminus \bigcup_{(a,b) \in \Gamma_q}
(a,b)+\CC\big)
 \bigg\},
\end{eqnarray*}
that is,
\[
 \ehk(A)  =
\lim_{q\to \infty} \frac{1}{q^3}\left[
  \sum_{d=0}^{\infty}
  \#\bigg\{\bbZ^2 \cap 
\bigg(d\PP \setminus \!\!\bigcup_{(a,b) \in \Gamma_q}
\!\!(a,b)
+\max\{0,d-q\}\PP\bigg)\bigg\}\right].
\]
\par
Also, we define a real continuous 
function $f : [0,\infty) \to \bbR$ by
\[
 f(t) = \;\text{the volume of }\; \left[t\PP \setminus
\!\!\bigcup_{(a,b) \in \Gamma} \!\!(a,b)
+\max\{0,t-1\}\PP \right]\;\text{in $\bbR^2$},
\]
then $\ehk(A) 
= \displaystyle{\int_0^{\infty}} f(t)\,dt$.
Let us denote the volume of $M \subseteq \bbR^2$ 
by $\vol(M)$.
To calculate $\ehk(A)$, 
we need to determine $f(t)$.
Namely, we need to show the following claim.

\begin{description}
\item[Claim]
\[
 f(t) =
\left\{\begin{array}{ll}
 \vol(t\PP) & (0 \le t < 1) \\[4pt]
 \vol(t\PP) - (n+4)\vol((t-1)\PP)
    & \left(1 \le t < \frac{n+2}{n+1}\right)\\[4pt]
 \frac{(n+2)t(2-t)}{2} + (n+2)\frac{(2-t)^2}{2n}
    & \left(\frac{n+2}{n+1} \le t < 2\right) \\[4pt]
 0 & (t \ge 2)
\end{array}
\right.
\]
\end{description}
\par
To prove the claim, we may assume that $t \ge 1$.
For simplicity, we put $M_{a,b} = (a,b)+ (t-1)\PP$ for every
point $(a,b)\in \Gamma$.
First suppose that $1 \le t < \frac{n+2}{n+1}$.
Then since $1 -n(t-1) > t-1$, 
$M_{0,0} \cap M_{1,0} = \emptyset$.
Similarly, one can easily see that any two $M_{a,b}$ 
do not intersect each
other.
Thus $f(t) = \vol(t\PP) - (n+4)\vol((t-1)\PP)$.
Next suppose that $\frac{n+2}{n+1} \le t < 2$.
Then $\PP \cap \{(x,y) \in \bbR^2\,|\, 0 \le y \le t-1\} 
= M_{0,0}\cup M_{1,0}\cup T_0$, 
where $T_0$ is a triangle with vertex
$(t-1,0)$, $(1,0)$ and $\left(t-1,\frac{2-t}{n}\right)$.
Similarly, there exist $(n+1)$ triangles 
$T_1,\ldots,T_{n+1}$
having the same volumes as $T_0$ such that
\[
\PP \cap \{(x,y)\in \bbR^2\,|\, 1 \le y \le t\}
= M_{-n,1} \cup \cdots \cup M_{1,1} 
\cup M_{0,1}\cup M_{1,1} \cup
T_1 \cup \cdots \cup T_{n+1}
\]
and any two $T_i$'s do not intersect each other.
Thus
\begin{eqnarray*}
f(t)
& = & 
\vol(\PP \cap \{(x,y)\in \bbR^2 \,|\, t-1 \le y \le 1\})
 +(n+2) \vol(T_0) \\
& = & \frac{(n+2)t(2-t)}{2} + (n+2)\frac{(2-t)^2}{2n}.
\end{eqnarray*}
Finally, suppose that $t \ge 2$. 
Then since $\PP$ is covered by $M_{a,b}$'s,
we have $f(t) =0$, as required.


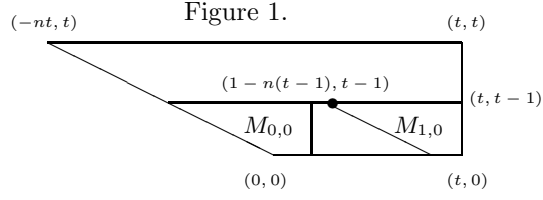
\begin{center}
\setlength{\unitlength}{1mm}
\begin{figure}[ht]
\par 
Figure 1. \\
\begin{picture}(100,22)
\put(80,5){\line(0,1){15}} 
\put(60,5){\line(0,1){7}} 
\put(55,5){\line(1,0){25}} 
\put(25,20){\line(1,0){55}} 
\put(41,12){\line(1,0){39}} 
\put(25,20){\line(2,-1){30}} 
\put(62,12){\line(2,-1){14}} 
\put(62,11){$\bullet$}
\put(48,14){\tiny $(1-n(t-1),t-1)$}
\put(20,22){\tiny $(-nt,t)$}
\put(78,22){\tiny $(t,t)$}
\put(81,12){\tiny $(t,t-1)$}
\put(51,1){\tiny $(0,0)$}
\put(78,1){\tiny $(t,0)$}
\put(51,8){\small $M_{0,0}$}
\put(71,8){\small $M_{1,0}$}
\end{picture}

\par 
The case where $1 \le t < \frac{n+2}{n+1}$ 
\end{figure}
\end{center}

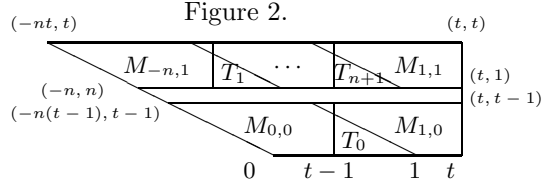
\begin{center}
\setlength{\unitlength}{1mm}
\begin{figure}[ht]
\par 
Figure 2. \\
\begin{picture}(100,22)
\put(80,5){\line(0,1){15}} 
\put(63,5){\line(0,1){7}} 
\put(63,14){\line(0,1){6}} 
\put(47,14){\line(0,1){6}} 
\put(55,5){\line(1,0){25}} 
\put(25,20){\line(1,0){55}} 
\put(37,14){\line(1,0){43}} 
\put(41,12){\line(1,0){39}} 
\put(25,20){\line(2,-1){30}} 
\put(44,20){\line(2,-1){12}} 
\put(60,12){\line(2,-1){14}} 
\put(60,20){\line(2,-1){12}} 
\put(20,22){\tiny $(-nt,t)$} 
\put(24,13){\tiny $(-n,n)$}
\put(20,10){\tiny $(-n(t-1),t-1)$}
\put(78,22){\tiny $(t,t)$}
\put(81,12){\tiny $(t,t-1)$}
\put(81,15){\tiny $(t,1)$}
\put(55,16){\small $\cdots$}
\put(51,2){\small $0$}
\put(59,2){\small $t-1$}
\put(73,2){\small $1$}
\put(78,2){\small $t$}
\put(64,6){\small $T_0$}
\put(48,15){\small $T_1$}
\put(63,15){\small $T_{n+1}$}
\put(51,8){\small $M_{0,0}$}
\put(71,8){\small $M_{1,0}$}
\put(35,16){\small $M_{-n,1}$}
\put(71,16){\small $M_{1,1}$}
\end{picture}
\par 
The case where $\frac{n+2}{n+1} \le t < 2$ 
\end{figure}
\end{center}

\par
Using the above claim, let us calculate $\ehk(A)$.
Note that $\vol(t\PP)= \frac{(n+2)t^2}{2}$.
\begin{eqnarray*}
 \ehk(A)
 & = &
 \int_{0}^{\frac{n+2}{n+1}} \frac{(n+2)t^2}{2}\, dt
 -(n+4) \int_{1}^{\frac{n+2}{n+1}} 
 \frac{(n+2)(t-1)^2}{2}\, dt \\
 & & + \int_{\frac{n+2}{n+1}}^{2} \frac{(n+2)t(2-t)}{2}\,dt
     + (n+2) \int_{\frac{n+2}{n+1}}^{2}
  \frac{(2-t)^2}{2n}\,dt \\
 & = & (n+2) \left[\frac{~1~}{2}+   \frac{1}{6(n+1)}\right],
\end{eqnarray*}
as required.
\end{proof}

\begin{discuss} \label{RNS-Poss}
Let $(A,\frm,k)$ be a complete unmixed local ring of
 positive characteristic
with $\dim A = 3$ and $e:=e(A) =3$.
What is the smallest value of $\ehk(A)$ among such rings?
\par
The function $f_e(s) = 3\left(\frac{s^3}{6}-
5 \frac{(s-1)^3}{6}\right)$,
which appeared in Eq.~(\ref{Ineq-gen}),
takes a maximal value
\[
f\left(\frac{5+\sqrt{5}}{4}\right) 
= \frac{15+5\sqrt{5}}{16} = 1.636\ldots
\]
in $s \in [1,2]$.
Hence $\ehk(A) \ge 1.636\ldots$.
But we beleive that this is not best possible.
\par
Suppose that $A$ is complete and $\ehk(A) <2$.
Then $A$ is F-rational by Lemma \ref{Mul34}.
Thus it is Cohen--Macaulay and 
$3+1 \le v =\emb(A) \le d+e-1 = 3+3-1=5$.
If $v \ne 5$, then $A$ is a hypersurface and
$\ehk(A) \ge \frac{~2~}{3}\cdot e = 2$ by 
Fact \ref{Hyp}.
Hence we may assume that $v =5$,
that is, $A$ has a maximal embedding dimension.
If we write as $A = R/I$, where $R$ is a complete 
regular local ring with
$\dim R = 5$, then $\height I =2$.
By Hilbert--Burch's theorem, there 
exists a $2\times 3$-matrix
$\bbM$ such that $I = I_2(\bbM)$,
the ideal generated by all $2$-minors of $\bbM$.
In particular, $A$ can be written as $A=B/aB$,
where $B = k[X]/I_2(X)$, $X$ is a generic $2\times 3$-matrix
and $a$ is a prime element of $B$.
This implies that
\[
 \ehk(A) = \ehk(B/aB) \ge \ehk(B)
= 3 \left\{\frac{1}{2}+\frac{1}{4!}\right\} = \frac{13}{~8~}=1.625;
\]
see \cite[Section 3]{EY}.
\par
For example, if $A= k[[T,xT,xyt,yT,x^{-1}yT]]$
is a rational normal scroll,
then $\ehk(A) = \frac{7}{~4~}=1.75$ by Example \ref{RNS}.
Is this the smallest value?
\end{discuss}

\begin{discuss} \label{Vero-Poss}
Let $(A,\frm,k)$ be a complete unmixed local ring of
positive characteristic.
Suppose that $\dim A =3$ and $e(A) =4$.
What is the smallest value of $\ehk(A)$ among such rings?
\par
As in Discussion \ref{RNS-Poss}, 
it suffices to consider F-rational local rings only.
For example, let $A = k[[x,y,z]]^{(2)}$ be 
the Veronese subring.
Then $A$ is an F-rational local
domain with $e(A) =4$  and $\ehk(A) =2$.
Also, let $A$ be the completion of the
Rees algebra $R(\frn)$ over
an F-rational double point $(R,\frn)$ of dimension $2$.
Then $A$ is an F-rational local domain
with $e(A) =4$ and $\ehk(A) \ge 2$
(we beleive that this inequality is strict).
\par
On the other hand,
the function $f_e(s)$ which appeared in Eq.~(\ref{Ineq-gen}),
takes a maximal value
\[
f\left(\frac{6+\sqrt{6}}{5}\right) 
= \frac{28+8\sqrt{6}}{25} = 1.903\ldots
\]
in $s \in [1,2]$.
Hence the fact that we can prove now is 
\lq\lq $\ehk(A) \ge 1.903\ldots$''
only.
\end{discuss}

\par
Based on Corollary \ref{2dimHK} and 
Discussion \ref{Vero-Poss}, we pose 
the following conjecture.

\begin{conj} \label{Vero-Conj}
Let $A$ be a complete unmixed local ring of 
positive characteristic with
$\dim A =3$. Let $r$ be an integer.
If $e(A) =r^2$, then
\[
  \ehk(A) \ge \frac{(r+1)(r+2)}{6}.
\]
Also, the equality holds if and only if $A$ is
isomorphic to $k[[x,y,z]]^{(r)}$.
\end{conj}

\par \vspace{2mm}
In the rest of this section, we prove the second statement of
Theorem \ref{Main3}.
Let $(A,\frm,k)$ be a complete unmixed local ring of 
characteristic $p>0$.
If $\ehk(A) =\frac{4}{~3~}$,
then $A$ is an F-rational hypersurface with $e(A) =2$
by the above observation.
Furthermore, suppose that $k = \overline{k}$ and 
$\chara k \ne 2$.
Then we may assume that $A$ can be written as the form
$k[[X,Y,Z,W]]/(X^2-\varphi(Y,Z,W))$.
To study Hilbert-Kunz multiplicities for these rings,
we prove the improved version of Theorem \ref{Key}.

\begin{prop} \label{Hyp2-Est}
Let $k$ be an algebraically closed field of $\chara k \ne 2$, 
and let $A = k[[X,Y,Z,W]]/(X^2-\varphi(Y,Z,W))$ be 
an F-rational hypersurface local ring.
Let $a$, $b$, $c$ be integers with $2 \le a \le b \le c$.
\par
Suppose that there exists a function $\ord: A \to \bbQ$ 
which satisfies the following conditions$:$
\begin{enumerate}
 \item $\ord(\alpha) \ge 0$; and
  $\ord(\alpha) = \infty \; \Longleftrightarrow \;\alpha = 0$.
 \item $\ord(x) = 1/2$, $\ord y = 1/a$, $\ord z = 1/b$, and
 $\ord w = 1/c$.
 \item $\ord(\varphi) \ge 1$.
 \item $\ord(\alpha+\beta)\ge 
   \min\{\ord(\alpha),\ord(\beta)\}$.
 \item $\ord(\alpha\beta) \ge \ord(\alpha) + \ord(\beta)$.
\end{enumerate}
Then we have
\[
  \ehk(A) \ge 2 - \frac{abc}{12}(N^3-n^3),
\]
where
\[
 N = \frac{~1~}{a} + \frac{~1~}{b}
 + \frac{~1~}{c} -\frac{~1~}{2},
 \qquad
 n = \max\left\{0,\;\;N- \frac{~2~}{c} \right\}.
\]
\par
In particular, if $(a,b,c) \ne (2,2,2)$, then $\ehk(A) > \frac{4}{~3~}$.
\end{prop}

\begin{proof}[\quad Proof]
First, we define a filtration $\{F_n\}_{n \in \bbQ}$ 
as follows$:$
\[
 F_{n} := \left\{\alpha \in A \,|\, 
  \ord(\alpha) \ge n \right\}.
\]
Then every $F_{n}$ is an ideal
and $F_{m}F_{n} \subseteq F_{m+n}$ holds 
for all $m,\,n \in \bbQ$.
Using $F_n$ instead of $\frm^n$, we shall estimate
$l_A(\frm^{[q]}/J^{[q]})$.
\par
Set $J = (y,z,w)A$ and fix a sufficiently large 
power $q=p^e$.
Put
\[
 s = \frac{1}{~a~}+\frac{1}{~b~}+\frac{1}{~c~},\qquad
 N = \frac{~1~}{a} + \frac{~1~}{b}
+ \frac{~1~}{c} -\frac{~1~}{2}.
\]
Since $J$ is a minimal reduction of $\frm$
and $xy^{q-1}z^{q-1}w^{q-1}$ generates the socle of 
$A/J^{[q]}$,
we have that $F_{sq} \subseteq J^{[q]}$.
Also, since $B=A/J^{[q]}$ is an Artinian 
Gorenstein local ring, we get
\[
 F_{\frac{(N+1)q}{2}}B \subseteq 0 :_B F_{\frac{Nq}{2}}B
 \cong K_{B/F_{\frac{Nq}{2}}B}.
\]
Hence, by Matlis duality theorem, we get
\[
 l_A \left(\frac{F_{\frac{(N+1)q}{2}}+J^{[q]}}{J^{[q]}}\right)
 \le l_B \left(K_{B/F_{\frac{Nq}{2}}B}\right)
 = l_B\left(B/F_{\frac{Nq}{2}}B\right).
\]
On the other hand, since $x^q \in F_{\frac{q}{2}}$ 
by the assumption, we have
\[
 x^q F_{\frac{Nq}{2}} \subseteq F_{\frac{(N+1)q}{2}}.
\]
Therefore by the similar argument 
as in the proof of Theorem \ref{Key}, we get

\begin{eqnarray*}
 l_A(\frm^{[q]}/J^{[q]})
& \le & l_A \left(
 \frac{Ax^q+J^{[q]}+F_{\frac{(N+1)q}{2}}}%
 {F_{\frac{(N+1)q}{2}}+J^{[q]}}\right)
 + l_A \left(
 \frac{F_{\frac{(N+1)q}{2}}+J^{[q]}}{J^{[q]}}\right) \\
& \le & l_A\left(A/(J^{[q]} +F_{\frac{(N+1)q}{2}}):x^q
 \right)
 + l_B\left(B/F_{\frac{Nq}{2}}B\right)  \\
& \le & 2 \cdot l_A\left(A/J^{[q]} + F_{\frac{N}{2}q} 
\right).
\end{eqnarray*}
\par
In fact, since
\begin{eqnarray*}
 & & \lim_{q \to \infty} \frac{1}{q^3} \;
l_A\left(A/J^{[q]}+F_{\frac{Nq}{2}}\right) \\
 & = & e(A) \cdot \lim_{q \to \infty} \frac{1}{q^3} \;\;
\vol \left\{(x,y,z) \in [0,q]^3 \,\bigg|\,
\frac{y}{~a~} + \frac{z}{~b~} + \frac{w}{~c~} 
\le \frac{Nq}{2} \right\} \\
 & = & 2 \cdot \vol \left\{(x,y,z) \in [0,1]^3 \,\bigg|\,
\frac{y}{~a~} + \frac{z}{~b~} + \frac{w}{~c~} 
\le \frac{N}{2} \right\} \\
 & = &   \frac{abc}{24}(N^3-n^3),
\end{eqnarray*}
we get
\[
 \ehk(A) \ge 2 - 2 \cdot \frac{abc}{24}(N^3-n^3) = 2 -
\frac{abc}{12}(N^3-n^3),
\]
as required.
\end{proof}

\begin{exam} \label{Special}
Let $k$ be an algebraically closed field 
of $\chara k \ne 2$.
Put $A = k[[X,Y,Z,W]]/(f(X,Y,Z,W))$.
\begin{eqnarray*}
 f(X,Y,Z,W)  = X^2+Y^3+Z^3+W^3
 & \Longrightarrow & \ehk(A) \ge \frac{55}{32};  \\[2mm]
 f(X,Y,Z,W)  = X^2+Y^2+Z^3+W^3
 & \Longrightarrow & \ehk(A) \ge \frac{14}{9};  \\[2mm]
 f(X,Y,Z,W)  = X^2+Y^2+Z^2+W^c
 & \Longrightarrow & \ehk(A) \ge \frac{9c^2-4}{6c^2}.
\end{eqnarray*}
\end{exam}

\par \vspace{2mm}
\par \noindent
\begin{proof}[\bf Proof of Theorem \ref{Main3}(2).]
Put $G = \gr_{\frm}(A)$ and $\fraM = \gr_{\frm}(A)_{+}$. 
The implication $(a) \Longrightarrow (b)$
follows from Proposition \ref{Hyp2-Est}.
$(b) \Longrightarrow (c)$ is clear.
Suppose $(c)$.
Then $\ehk(G_{\fraM}) = \frac{4}{~3~}$.
Also, by Proposition \ref{grHK} and Theorem \ref{Main3} (1),
we have that $ \frac{4}{~3~} \le \ehk(A) \le \ehk(G_{\fraM}) 
=  \frac{4}{~3~}$.
Thus $\ehk(A) = \frac{4}{~3~}$, as required.
\end{proof}

\par
Also, the following corollary follows from the proof of Proposition 
\ref{Hyp2-Est} and Example \ref{Special}.

\begin{cor} \label{A_1}
Let $(A,\frm,k)$ be an unmixed local ring of
characteristic $p >0$ with $\dim A =3$.
Assume that $k = \overline{k}$ and $p \ne 2$.
Then the following conditions are equivalent$:$
\begin{enumerate}
 \item $\frac{4}{3} < \ehk(A) \le \frac{3}{2}$.
 \item $\gr_{\frm}(A) \cong k[X,Y,Z]/(X^2+Y^2+Z^2)$.
 \item $A$ is isomorphic to a hypersurface 
   $k[[X,Y,Z,W]]/(X^2+Y^2+Z^2+W^c)$
  for some integer $c \ge 3$.
\end{enumerate}
When this is the case, $\ehk(A) \ge \dfrac{~3~}{2} 
- \dfrac{~2~}{3c^2}$.
\end{cor}
\section{A generalization of the main result to higher 
dimensional case}

\par
In this section, we want to consider a generalization of 
Theorem \ref{Main3} in case of $\dim A \ge 4$.
Let $d \ge 1$ be an integer and $p> 2$ a prime number.
If we put
\[
 A_{p,d} := 
\overline{\bbF_p}[[X_0,X_1,\ldots,X_d]]/(X_0^2+\cdots + X_d^2),
\]
then we guess that $\ehk(A_{p,d}) = \shk(p,d)$ 
holds according to the
observations until the previous section.
In the following, let us formulate this as a conjecture 
and prove that it is also true in case of $\dim A =4$.
\par
In \cite{HM}, Han and Monsky gave an algorism 
to calculate $\ehk(A_{p,d})$,
but it is not so easy to represent it as a quotient of 
two polynomials of $p$ for any fixed $d \ge 1$.
\[
\begin{array}{|c||c|c|c|c|} \hline
d & 1 & 2 & 3 & 4 \\ \hline
\ehk(A_{p,d}) & 2 & \frac{~3~}{2} & \frac{~4~}{3}
& \frac{29p^2+15}{24p^2+12} \\ \hline
\end{array}
\]
\par
On the other hand, surprisingly, 
Monsky proved the following theorem$:$

\begin{thm}[Monsky {\rm \cite{Mo2}}]
Under the above notation, we have
\begin{equation}
 \lim_{p \to \infty} \ehk(A_{p,d})
 = 1 + \frac{c_d}{d!},
\end{equation}
where
\begin{equation}
 \sec x + \tan x = \sum_{d=0}^{\infty} \frac{c_d}{d!}x^d
 \quad \left(|x| < \frac{\pi}{~2~} \right). 
\end{equation}
\end{thm}

\begin{remark} \label{SecTan}
It is known that the Taylor expansion of $\sec x$ 
(resp. $\tan x$) at origin
can be written as follows$:$
\begin{eqnarray*}
 \sec x & = & \sum_{i=0}^{\infty} 
  \frac{E_{2i}}{(2i)!} x^{2i}, \\
 \tan x & = & \sum_{i=1}^{\infty} (-1)^{i-1}
\frac{2^{2i}(2^{2i}-1)B_{2i}}{(2i)!}x^{2i-1},
\end{eqnarray*}
where $E_{2i}$ (resp. $B_{2i}$) is said to be Euler number
(resp. Bernoulli number).
\par
Also, $c_d$ appeared in Eq.(4.1) is a positive integer 
since $\cos t$ is an unit element in a ring
$\HH = \left\{\sum_{n=0}^{\infty} a_n \frac{t^n}{n!}
\,|\, a_n \in \bbZ \;
\text{for all $n \ge 0$} \right\}$.
\end{remark}

\par
Based on the above observation, we pose the 
following conjecture.

\begin{conj} \label{Conj-Gen}
Let $d \ge 1$ be an integer and $p > 2$ a prime number.
Put
\[
  A_{p,d} := \overline{\bbF_p}[[X_0,X_1,\ldots,X_d]]/
(X_0^2+\cdots + X_d^2).
\]
Let $(A,\frm,k)$ be a $d$-dimensional unmixed local ring
with $k = \overline{\bbF_p}$.
Then the following statements hold.
\begin{enumerate}
 \item If $A$ is not regular, then $\ehk(A) 
\ge \ehk(A_{p,d}) \ge 1+\frac{c_d}{d!}$.
In particular, $\shk(p,d) = \ehk(A_{p,d})$.
 \item If $\ehk(A) = \ehk(A_{p,d})$, 
then $\widehat{A} \cong A_{p,d}$
as local rings.
\end{enumerate}
\end{conj}

\par \vspace{2mm}
In the following, we prove that this is true 
in case of $\dim A =4$.
Note that
\[
 \lim_{p\to \infty} \ehk(A_{p,4})
= \lim_{p \to \infty} \frac{29p^2+15}{24p^2+12}
= \frac{29}{24}= 1 + \frac{c_4}{4!}.
\]

\begin{thm} \label{4Main}
Let $(A,\frm,k)$ be an unmixed local ring of 
characteristic $p>0$ with
$\dim A =4$.
If $e(A) \ge 3$, then $\ehk(A) \ge \frac{~5~}{4} 
= \frac{30}{24}$.
\par \vspace{2mm}
Suppose that $k = \overline{k}$ and $\chara k \ne 2$.
Put
\[
 A_{p,4} = \overline{\bbF_p}[[X_0,X_1,\ldots,X_4]]
/(X_0^2+\cdots + X_4^2).
\]
Then the following statement holds.
\begin{enumerate}
 \item If $A$ is not regular, then
\[
  \ehk(A) \ge \ehk(A_{p,4}) = \frac{29p^2+15}{24p^2+12}.
\]
 \item The following conditions are equivalent$:$
\begin{enumerate}
   \item Equality holds in $(1)$.
   \item $\ehk(A) < \frac{~5~}{4}$.
   \item The completion of $A$ is isomorphic to $A_{p,4}$.
\end{enumerate}
\end{enumerate}
\end{thm}

\begin{proof}[\quad Proof]
Put $e = e(A)$, the multiplicity of $A$.
We may assume that $A$ is complete with $e \ge 2$ and 
$k$ is infinite.
In particular, $A$ is a homomorphic 
image of a Cohen--Macalay local ring, and 
there exists a minimal reduction $J$ of $\frm$.
Then $\mu_A(\frm/J^{*}) \le e-1$ by 
Lemma \ref{Sally-type}.
We first show that $\ehk(A) \ge \frac{5}{~4~}$ 
if $e \ge 3$.

\vspace{2mm}
\begin{description}
 \item[Claim 1] If $3 \le e \le 10$, 
then $\ehk(A) \ge \frac{5}{~4~}$.
\end{description}

Putting $r = e-1$ and $s = 2$ in Theorem \ref{Key}, 
since $v_2 = \frac{1}{~2~}$,
we have
\[
\ehk(A) \ge e \left\{v_2 - \frac{(e-1)1^4}{4!}\right\}
= \frac{(13-e)e}{24} \ge \frac{30}{24},
\]
as required.

\vspace{2mm}
\begin{description}
 \item[Claim 2] If $11 \le e \le 29$, 
then $\ehk(A) \ge \frac{737}{384}
\left(> \frac{5}{~4~}\right)$.
\end{description}

By Fact \ref{Hyp}, we have 
$v_{3/2} = \frac{1-\beta_{4+1}}{2} =
\frac{77}{384}$.
Putting $r=e-1$ and $s = \frac{3}{2}$ in Theorem \ref{Key},
we have
\[
  \ehk(A)
 \ge e \left\{v_{3/2} - \frac{e-1}{24} \cdot 
\left(\frac{1}{2}\right)^4
\right\}
 = \frac{(78-e)e}{384} \ge \frac{11(78-11)}{384} 
= \frac{737}{384},
\]
as required.

\vspace{2mm}
\begin{description}
 \item[Claim 3] If $e \ge 30$, then $\ehk(A) 
\ge \frac{5}{~4~}$.
\end{description}

By Proposition \ref{IneqMul}, we have $\ehk(A) 
\ge \frac{e}{4!} \ge
\frac{30}{24}$.

\par \vspace{2mm}
In the following, we assume that $k = \overline{k}$, 
$\chara k \ne 2$ and $e
\ge 2$.
To see (1),(2), we may assume that $e = 2$ by the 
above argument.
Then since $\ehk(A) = 2$ if $A$ is not F-rational, 
we may also assume that
$A$ is
F-rational and thus is a hypersurface.
Thus $A$ can be written as the following form$:$
\[
  A = k[[X_0,X_1,\ldots,X_4]]/
(X_0^2- \varphi(X_1,X_2,X_3,X_4))
\]
If $A$ is isomorphic to $A_{p,4}$, then by \cite{HM}, 
it is known that
\[
 \ehk(A) = \frac{29p^2+15}{24p^2+12}.
\]
Suppose that $A$ is not isomorphic to $A_{p,4}$.
Then one can take a minimal numbers of generators 
$x,y,z,w,u$ of $\frm$ and
one can define a function $\ord :A \to \bbQ$ such that
\[
 \ord(x) =  \ord(y) =  \ord(z) =  \ord(z) =  \frac{1}{~2~},
\quad  \ord(u) = \frac{1}{~3~}.
\]
If we put $J = (y,z,w,u)A$ and
$F_n = \{\alpha \in A \,|\, \ord(\alpha) \ge n\}$, then
by the similar argument as in the proof of Proposition 
\ref{Hyp2-Est},
we have
\[
  l_A(\frm^{[q]}/J^{[q]}) \le 2  \cdot l_A(A/J^{[q]}+F_{2q/3}).
\]
Divided the both-side by $q^d$ and 
taking a limit $q \to \infty$, we get
\[
\e(A) - \ehk(A) \le 2 \cdot 
e(A)\cdot \vol\left\{(y,z,w,u) \in [0,1]^4
\,\bigg|\,
  \frac{y}{2}+\frac{z}{2}+\frac{w}{2}+\frac{u}{3} 
\le \frac{2}{3}\right\}.
\]
To calculate the volume in the right-hand side, we put
\[
 F_u =
 \left\{
 \begin{array}{ll}
 \frac{1}{6}\left(\frac{4}{3}-\frac{2}{3}u\right)^3 -
         \frac{1}{6}\left(\frac{1}{3}-\frac{2}{3}u\right)^3
     & (0 \le u \le \frac{1}{2}) \cr
      \frac{1}{6}\left(\frac{4}{3}-\frac{2}{3}u\right)^3
     & (\frac{1}{2} \le u \le 1)
 \end{array}
\right.
\]
Then one can easily calculate
\[
 \text{the above volume} =  \int_{0}^{1} F_u \,du 
  = \frac{237}{2^4 3^4}.
\]
It follows that
\[
 \ehk(A) \ge 2 - 4 \times \frac{237}{2^4 3^4} 
 = \frac{411}{324} >
\frac{5}{~4~}.
\]
\end{proof}

\par
The following conjecture also holds if $\dim A \le 4$.

\begin{conj} \label{Conj-add}
Under the same notation as in Conjecture $\ref{Conj-Gen}$,
if $e(A) \ge 3$, then
\[
  \ehk(A) \ge \frac{c_d+1}{d!}.
\]
\end{conj}

\par
\begin{discuss} \label{Homo-Conj}
Let $d \ge 2$ be an integer and fix a prime number 
$p \gg d$.
Assume that Conjectures \ref{Conj-Gen} and \ref{Conj-add} 
are true.
Also, assume that $\shk(p,d) < \shk(p,d-1)$ for all 
$d \ge 3$.
Let $A = k[X_0,\ldots,X_v]/I$ be a $d$-dimensional
homogeneous unmixed $k$-algebra with $\deg X_i =1$,
and let $\frm$ be a unique homogeneous maximal 
ideal of $A$. 
Suppose that $\chara k =p > 2$ and 
$k=\overline{\bbF_p}$.
Then $\ehk(A) = \shk(p,d)$ implies that
$\widehat{A_{\frm}} \cong A_{p,d}$.
\par
In fact, if $\ehk(A) = \shk(p,d)$, then we may assume that
$\ehk(A) \le \frac{c_d+1}{d!}$.
Thus $e(A_{\frm}) =2$ by the assumption that 
Conjecture \ref{Conj-add}
is true.
For any prime ideal $PA_{\frm}$ of $A_{\frm}$ such that
$P \ne \frm$, we have 
$\ehk(A_P) \le \ehk(A) = \shk(p,d) < \ehk(p,n)$, where
$n = \dim A_P <d$.
Since $A_P$ is also unmixed,
we have that $A_P$ is regular.
Thus $A_{\frm}$ has an isolated singularity.
Hence $A$ is a non-degenerate quadric hyperplane.
In other words, $\widehat{A_{\frm}}$ is isomorphic to $A_{p,d}$.
\end{discuss}
\bibliographystyle{amsplain}

\end{document}